\documentclass[a4paper,12pt]{amsart}
\usepackage{amsmath}
\usepackage{amssymb}
\usepackage{amsthm}
\usepackage{graphicx}

\newtheorem{theorem}{\indent Theorem}
\newtheorem{lemma}{\indent Lemma}
\newtheorem{proposition}{\indent Proposition}

\newtheorem{definition}{\indent Definition}

\newtheorem{example}{Example}

\providecommand{\abs}[1]{\lvert#1\rvert}
\providecommand{\norm}[1]{\displaystyle\left\lVert#1\right\rVert}
\providecommand{\keywords}[1]{\textbf{\textit{Keywords:}} #1.}

\begin{document}

\title{Averaging method for dynamic systems \\ on time scales with periodicity}
\author{Aleksey Ogulenko}
\email{ogulenko.a.p@onu.edu.ua}
\address{I. I. Mechnikov Odesa National University, Dvoryanska str, 2, Odesa, Ukraine, 65082}

\begin{abstract}
	This paper aims to improve existing results about using averaging
	method for analysis of dynamic systems on time scales.  We obtain a
	more accurate estimate for proximity between solutions of original and
	averaged systems regarding $\Delta$--periodic and
	$\Delta$-quasiperiodic systems, which are introduced for the first
	time. To illustrate the application of the averaging theorem for such
	kind of system we considered an example and conducted numerical
	modelling.  Obtained results extend an application area for previously
	developed numer\-ical\-ly--asymptotic method of solution for optimal
	control problems on time scales.
\end{abstract}


\maketitle

\section{Introduction}
A systematic theory of averaging method for ordinary differential equations 
began from the works of \cite{Krylov}.
Further it was developing by \cite{Bogolubov} and the others. 
Since then, there have been many works establishing the averaging method 
for various types of dynamic systems: differential equations with discontinuous and multi-valued 
right-hand side, with Hukuhara derivative, with delay etc.
The review of these results one can find in \cite{Plotnikov}.

On the other hand, the theory of dynamic equations on time scales was
introduced by \cite{Hilger1988} in order to unify continuous and discrete
calculus.  In detail, the description of time scale analysis can be found in
the \cite{BohnerPeterson2001, BohnerPetersonAdvances}.

As far as we know the averaging method in connection with the systems on time scales
was first examined by \cite{SlavikAveragin}. 
In particular, there were studied conditions of proximity between solutions of the
original system on time scale and some generalized differential equation.
From the practical point of view, the interpretation of the last equation's 
solution in terms of given application is somewhat unclear.

Previously we established the scheme of full averaging for dynamic systems on
time scales (\cite{Ogulenko2012}) and in our approach the averaged system has
the same time nature as the original one.  Recently we also established the
analogous result for partially averaged systems, \cite{OgulenkoPartial2017}. On
the base of this scheme, it was developed the numer\-ical\-ly--asymptotic
method of solution for optimal control problems on time scales
(\cite{OgulenkoKichmarenko2016,OgulenkoKichmarenko2017}).  

\section{Preliminary results.}

We now present some basic information about time scales according to
\cite{BohnerPeterson2001}.  A time scale is defined as a nonempty closed subset
of the set of real numbers and usually denoted by $\mathbb{T}$.  The properties
of the time scale are determined by the following three functions:
1) the forward-jump operator $\sigma(t) = \inf \left\{s \in \mathbb{T}: s > t\right\}$;
2) the backward-jump operator $\rho(t) = \sup \left\{s \in \mathbb{T}: s < t \right\}$ 
(in this case, we set $\inf\varnothing  = \sup\mathbb{T}$ and $\sup\varnothing = \inf\mathbb{T}$);
3) the granularity function $\mu(t) = \sigma(t) - t$.

The behaviour of the forward- and backward-jump operators at a given point of
the time scale specifies the type of this point. If $t < \sigma(t)$, then $t$
is called right--scattered, if $t = \sigma(t)$ --- right--dense. Also, point
will be called left--scattered when $\rho(t) < t$ and  left--dense when
$\rho(t) = t$. Finally, point is called dense if it is right-dense and
left--dense at the same time and isolated if it is both right--scattered and
left-scattered.

Important role in time scales calculus has the set $\mathbb{T}^\kappa$ which is
derived from the time scale $\mathbb{T}$ as follows: if $\mathbb{T}$ has a
left--scattered maximum $m$, then $\mathbb{T}^\kappa = \mathbb{T} - \{m\}$.
Otherwise, $\mathbb{T}^\kappa = \mathbb{T}$.  In what follows, we set $\left[a,
b\right]_\mathbb{T} = \left\{t \in \mathbb{T}: a \leqslant t \leqslant b
\right\}$.

\begin{definition}[\cite{BohnerPeterson2001}]
    Let $f:\mathbb{T} \rightarrow \mathbb{R}$ and $t \in \mathbb{T}^\kappa$.
    The number $f^{\Delta}(t)$ is called $\Delta$-derivative of function $f$ at the point $t$,
    if $\forall \varepsilon > 0$ there exists a neighborhood $U$ of the point $t$ 
    (i.~e., $U = (t - \delta, t + \delta) \cap \mathbb{T}, \delta < 0$) such that 
    \[
    \abs{f(\sigma(t)) - f(s) - f^{\Delta}(t)(\sigma(t)-s)} \leqslant
        \varepsilon\abs{\sigma(t) - s} \quad \forall s \in U.
    \]
\end{definition}

\begin{definition}[\cite{BohnerPeterson2001}]
    If $f^{\Delta}(t)$ exists $\forall t \in \mathbb{T}^\kappa$,
    then $f:\mathbb{T} \rightarrow \mathbb{R}$ is called $\Delta$-differentiable on 
    $\mathbb{T}^\kappa$. The function $f^{\Delta}(t): \mathbb{T}^\kappa \rightarrow \mathbb{R}$
    is called the delta-derivative of a function $f$ on  $\mathbb{T}^\kappa$. 
\end{definition}
If $f$ is differentiable with respect to $t$, then $f(\sigma(t)) = f(t) + \mu(t)f^{\Delta}(t)$.

\begin{definition}[\cite{BohnerPeterson2001}]
	The function $f:\mathbb{T} \rightarrow \mathbb{R}$ is called
	$rd$-continuous if it is continuous at the right-dense points and has
	finite left limits at the left-dense points.  The set of these
	functions is denoted by $C_{rd} = C_{rd}(\mathbb{T}) = 
	C_{rd}(\mathbb{T}; \mathbb{R})$.
\end{definition}

The indefinite integral on the time scale takes the form $\int{f(t) \Delta t} =
F(t) + C$, where $C$ is integration constant and $F(t)$ is the preprimitive for
$f(t)$.  If the relation $F^{\Delta}(t) = f(t)$ where
$f:\mathbb{T}\rightarrow\mathbb{R}$ is an $rd$-continuous function, is true for
all $t \in \mathbb{T}^\kappa$ then $F(t)$ is called the primitive of the
function $f(t)$. If $t_0 \in \mathbb{T}$ then $F(t) = \int\limits_{t_0}^t{f(s)
\Delta s}$ for all $t$.  The definite $\Delta$-integral on time scale interval
is defined by Newton--Leibniz formula.

\begin{definition}[\cite{BohnerPeterson2001}]
A function $p:\mathbb{T} \rightarrow \mathbb{R}$ is called regressive (positive regressive) if
    \[
        1 + \mu(t)p(t) \neq 0, \quad (1 + \mu(t)p(t) > 0), \qquad t \in \mathbb{T}^\kappa.
    \]
The set of regressive (positive regressive) and $rd$-continuous functions is denoted by 
$\mathcal{R} = \mathcal{R}(\mathbb{T})$ ($\mathcal{R}^+ = \mathcal{R}^+(\mathbb{T})$).
\end{definition}

A function $p$ from the class $\mathcal{R}$ can be associated
with a function $e_p(t, t_0)$ which is the unique solution of Cauchy problem
\[
    y^{\Delta} = p(t)y, \quad y(t_0) = 1.
\]
The function $e_p(t,t_0)$ is an analog, by its properties, of the exponential
function defined on~$\mathbb{R}$.

In what follows we heavily use the next result.
\begin{theorem}[Substitution rule, \cite{BohnerPeterson2001}, theorem 1.98]
	Assume $\nu : \mathbb{T} \rightarrow \mathbb{R}$ is strictly
	increasing and $\tilde{\mathbb{T}} = \nu(\mathbb{T})$ is a time scale.
	If $f: \mathbb{T} \rightarrow \mathbb{R}$ is an rd-continuous function and
	$\nu$ is differentiable with rd-continuous derivative, then for $a, b \in \mathbb{T}$,
	\begin{equation}
		\label{substitution-rule}
		\int\limits_a^b g(s)\nu^\Delta(s)\,\Delta{s} = 
		\int\limits_{\nu(a)}^{\nu(b)} g\left(\nu^{-1}(s)\right)\,\tilde\Delta{s}. 
	\end{equation}
\end{theorem}

Various kinds of periodicity on time scale was presented and studied by
\cite{Adivar2013}. The basic framework is as follows. For arbitrary non-empty
susbet $\mathbb{T}^*$ of the time scale $\mathbb{T}$ including a fixed number
$t_0$ the operators $\delta_\pm : [t_0, +\infty) \times \mathbb{T}^*
\rightarrow \mathbb{T}^*$ are introduced.  The operators $\delta_+$ and
$\delta_-$ associated with the initial point $t_0 \in \mathbb{T}^*$ are
said to be forward and backward shift operators on the set
$\mathbb{T}^*$, respectively. The first argument in $\delta_\pm(s, t)$
is called the shift size. The values $\delta_+(s, t)$ and
$\delta_-(s,t)$ indicate translation of the point $t \in \mathbb{T}^*$
to the right and left by $s$ units, respectively.

\begin{definition}[\cite{Adivar2013}]
	Let $\mathbb{T}$ be a time scale with the shift operators
	$\delta_\pm$ associated with the initial point $t_0 \in \mathbb{T}^*$.
	The time scale $\mathbb{T}$ is said to be periodic in shifts $\delta\pm$
	if there exists a $p \in (t_0, \infty)_{\mathbb{T}^*}$ such that 
	$(p, t) \in D_\mp$ for all $t \in \mathbb{T}^*$. Furthermore, if 
	\[
		P = \inf\left\{
			p \in (t_0, \infty)_{\mathbb{T}^*}: 
			(p, t) \in D_\mp \text{ for all } t \in \mathbb{T}^*
		     \right\} \neq t_0,
	\]
	then $P$ is called the period of the time scale $\mathbb{T}$.
\end{definition}

\begin{definition}[\cite{Adivar2013}]
	Let $\mathbb{T}$ be a time scale that is periodic in shifts $\delta_\pm$
	with the period $P$. We say that a real valued function 
	$f$ defined on $\mathbb{T}^*$ is periodic in shifts $\delta_\pm$ 
	if there exists a $T \in [P, \infty)_{\mathbb{T}^*}$ such that 
	$(T, t) \in D_\pm$ and $f(\delta_\pm(T, t)) = f(t)$ for all $t \in \mathbb{T}^*$.
	The smallest such a number $T \in [P, \infty)_{\mathbb{T}^*}$ 
	is called the period of~$f$.
\end{definition}

\begin{definition}[\cite{Adivar2013}]
	Let $\mathbb{T}$ be a time scale that is periodic in shifts $\delta_\pm$
	with the period $P$. We say that a real valued function 
	$f$ defined on $\mathbb{T}^*$ is $\Delta$-periodic in shifts $\delta_\pm$ 
	if there exists a $T \in [P, \infty)_{\mathbb{T}^*}$ such that 
	$(T, t) \in D_\pm$ for all $t \in \mathbb{T}^*$, the shifts 
	$\delta_\pm(T, t)$ are $\Delta$-differentiable with rd-continuous 
	derivative with respect to second argument and 
	\[
		f(\delta_\pm(T, t))\delta_\pm^\Delta(T, t) = f(t)
	\]
	for all $t \in \mathbb{T}^*$ The smallest such a number 
	$T \in [P, \infty)_{\mathbb{T}^*}$ is called the period of~$f$.
\end{definition}

It was shown in \cite{Adivar2013} that the following propositions about 
periodicity in shifts are true.

\begin{proposition}[\cite{Adivar2013}]
	If $\delta_+(s, \cdot)$ is $\Delta$-differentiable in its
	second argument, then $\delta_+^\Delta(s, \cdot) > 0$.
\end{proposition}

\begin{proposition}[\cite{Adivar2013}]
	Let $\mathbb{T}$ be a time scale that is periodic in shifts $\delta_\pm$
	with the period $P$ and $f$ a $\Delta$-periodic in shifts $\delta_\pm$ 
	with the period $T \in [P, \infty)_{\mathbb{T}^*}$.
	Suppose that $f \in C_{rd}(\mathbb{T})$, then
	\[
		\int\limits_{t_0}^{t} f(s) \Delta{s} = 
		\int\limits_{\delta^{T}_\pm(t_0)}^{\delta^{T}_\pm(t)} f(s) \Delta{s}.
	\]
\end{proposition}

\section{Main results}
	
Let $\mathbb{T}$ be an unbounded above time scale that is periodic in shifts
$\delta_\pm$ with period $P \in \left(t_0, +\infty\right)_{\mathbb{T}^*}$.  For
simplicity we denote by $\delta_{\pm}^T(t)$ the shift operators with
period $T$ and by $\delta^{(i)}_{\pm T}(t)$ or $\delta^{(i)}(t)$ the $i$-th
power of shift operator composition, dropping argument sometimes.

Consider on $\mathbb{T}$ the following dynamic system:
\begin{equation}
    \label{mainsys-periodic}
        x^{\Delta} = \varepsilon X(t, x), \quad x(t_0) = x_0.
\end{equation}
Here $x \in \mathbb{R}^n$, $\varepsilon > 0 $ is a small parameter, 
$X(t,x)$ is $n$-dimensional vector--function such that every component is 
$\Delta$-periodic in shifts $\delta_\pm(T, t)$ function, 
$T \in \left[P, +\infty\right)_{\mathbb{T}^*}$.

In correspondence to this original system, we put another dynamic system on the
same time scale as follows:
\begin{equation}
    \label{avgsys-periodic}
    \xi^{\Delta} = \varepsilon \widetilde{X}(t, \xi), \quad \xi(t_0) = x_0,
\end{equation}
where
\begin{equation}
    \label{avglimit-periodic}
    \begin{aligned}
    \widetilde{X}(t, x) = 
    	\Bigg\{
	\widetilde{X}_i(x) 
		&= \frac{1}{\delta^{(i+1)}(t_0) - \delta^{(i)}(t_0)} 
			\int\limits_{\delta^{(i)}(t_0)}^{\delta^{(i+1)}(t_0)} X(t,x) \Delta{t}, \\
		& \delta^{(i)}(t_0) \leqslant t < \delta^{(i+1)}(t_0), 
			\quad i = 0, 1, 2, \dots 
	\Bigg\}.
    \end{aligned}
\end{equation}
The last system \eqref{avglimit-periodic} we call partially averaged system corresponding to the original one.

Taking into account $\Delta$-periodical properties of the function $X(t, x)$,
it is easy to see, that
    \[
	    \widetilde{X}_i(\xi) = 
	    	\frac{\int\limits_{\delta^{(i)}(t_0)}^{\delta^{(i+1)}(t_0)} X(t,\xi) \Delta{t}}
			{\delta^{(i+1)}(t_0) - \delta^{(i)}(t_0)} = 
		\frac{\int\limits_{\delta^{(i-1)}(t_0)}^{\delta^{(i)}(t_0)} X(t,\xi) \Delta{t}}
			{\delta^{(i+1)}(t_0) - \delta^{(i)}(t_0)} = 
		\dots = 
		\frac{\int\limits^{\delta^{T}_{+}(t_0)}_{t_0} X(t,\xi) \Delta{t}}
			{\delta^{(i+1)}(t_0) - \delta^{(i)}(t_0)},
    \]
    that is, 
    \[
	    \widetilde{X}_i(\xi) = \frac{\delta^T_+(t_0) - t_0}{\delta^{(i+1)}(t_0) - \delta^{(i)}(t_0)}\widetilde{X}_0(\xi),
	    	\qquad i = 0, 1, 2, \dots.
    \]

We now prove that under general conditions there exists proximity between 
solutions of systems \eqref{mainsys-periodic} and \eqref{avgsys-periodic}.

\begin{theorem}
\label{theorem-periodic-averaging}
Let $Q = \Big\{t \in \mathbb{T}, x \in D \Big\}$, $x(t)$ and $\xi(t)$ denote 
solutions of the Cauchy problems \eqref{mainsys-periodic} and \eqref{avgsys-periodic} respectively. 
Now suppose the following conditions hold in $Q$:
\begin{itemize}
	\item[1)] every component of vector--function $X(t, x)$
		is $\Delta$-periodic in shifts $\delta_\pm(T, t)$ function, $T
		\in \left[P, +\infty\right)_{\mathbb{T}^*}$.
	\item[2)] the function $X(t,x)$ is rd-continuous with respect to $t$
		and regressive. Moreover, $X(t, x)$ satisfies conditions of
		existence and uniqueness of solution for Cauchy problem such
		that
		\[
			\forall (t,x) \in Q \quad \left\|X(t,x)\right\| \leqslant M, M~>~0,
		\]
		$X(t,x)$ is Lipschitz continuous with respect to $x$ with constant $\lambda~>~0$,~i.~e.
        	\[
	            \left\|X\left(t,x_1\right) - X\left(t,x_2\right)\right\| \leqslant 
		    	\lambda \left\|x_1 - x_2\right\| \qquad
        	        \forall \left(t,x_1\right), \left(t, x_2\right) \in Q ~;
	        \]
	\item[3)] there exists a constant $K > 0$ such that the following 
		holds for all $i \geqslant 1$:
		\[
			\delta^{(i+1)}(t_0) - \delta^{(i)}(t_0) \leqslant K;
		\]

        \item[4)] the solution $\xi(t)$ of averaged system \eqref{avgsys-periodic} with initial value 
		$\xi(t_0) = x_0 \in D' \subset D$ is well defined for all $t \in \mathbb{T}^\kappa$ 
		and with its $\rho$-neighbourhood lies in~$D$.
\end{itemize}
    Then for any $L > 0$ there exists $\varepsilon_0\left(L\right) > 0$ such that
    for $0 < \varepsilon <~\varepsilon_0$ and $t \in \left[t_0, t_0 + L\varepsilon^{-1}\right] \cap \mathbb{T}$ 
    the following estimate holds:
    \begin{equation}
    	\label{solution-proximity}
        \norm{x(t) - \xi(t)} \leqslant C\varepsilon.
    \end{equation}
\end{theorem}

\proof
    It is easy to see that $\widetilde{X}(t, x)$ is bounded and Lipschitz 
    continuous with respect to the second argument. It directly follows from the way of construction 
    \eqref{avglimit-periodic}. So, we have for any fixed $t \in \mathbb{T}$
    \[
        \norm{\widetilde{X}(t, x') - \widetilde{X}(t, x'')} \leqslant \lambda \norm{x' - x''}.
    \]
    Therefore, conditions 1) and 2) imply the existence and uniqueness of
    solutions for both original system and averaged one. Moreover, these
    solutions can be continued until $x(t) \in D$ (accordingly, $\xi(t) \in D$).

    Let us write both original and partially averaged systems in integral form:
    \begin{equation*}
        x(t) = x_0 + \varepsilon\int\limits_{t_0}^t{X(s, x(s))\Delta{s}}, \qquad
        \xi(t) = x_0 + \varepsilon\int\limits_{t_0}^t{\widetilde{X}(s, \xi(s))\Delta{s}}.
    \end{equation*}
    In the same way as we did establishing the scheme of full averaging for
    dynamic systems on time scales in \cite{Ogulenko2012}, let us estimate the
    norm of difference between solutions:
    \begin{eqnarray*}
        \norm{x(t) - \xi(t)} &=& \norm{\varepsilon \int\limits_{t_0}^t
		\left[X(s, x(s)) - \widetilde{X}(s, \xi(s))\right]\Delta{s}} \leqslant \\
        & \leqslant & \lambda\varepsilon\int\limits_{t_0}^t\norm{ x(s) -  \xi(s)}\Delta{s} +
                        \varepsilon\norm{\int\limits_{t_0}^t\left[X(s, \xi(s)) - \widetilde{X}(s, \xi(s))\right]\Delta{s} }~.
    \end{eqnarray*}
    We will estimate the last summand on the time scale interval 
    $\left[t_0, t_0 + L\varepsilon^{-1}\right]~\cap~\mathbb{T}$. 

    By $\varphi(t, \xi)$ denote the last integrand:
    \[
	    \varphi(t, \xi) = X(t, \xi(s)) - \widetilde{X}(t, \xi(s)).
    \]

    Consider time interval $\left[\delta^{(i)}(t_0), \delta^{(i+1)}(t_0)\right]$.
    By construction, on this interval $\widetilde{X}(t,\xi) = \widetilde{X}_i(\xi)$
    and 
    \[
	    \int\limits_{\delta^{(i)}(t_0)}^{\delta^{(i+1)}(t_0)} 
	    	\!\!\!\!\!\!\varphi(s, \xi_i) \Delta{s} = 0,
    \]
    where $\xi_i = \xi\left(\delta^{(i)}(t_0)\right) = const$.  
    
    Further, 
    \[
	\begin{gathered}
		\left\| 
			\int\limits_{t_0}^t \varphi(s, \xi(s)) \Delta{s} 
		\right\| \leqslant
		\left\| 
			\int\limits_{t_0}^{\delta_+^{(N)}(t_0)}	
			\!\!\!\!\!\! \varphi(s, \xi(s)) \Delta{s} 
		\right\| + 
		\left\|\ \int\limits_{\delta^{(N)}(t_0)}^t	
			\!\!\!\!\! \varphi(s, \xi(s)) \Delta{s} 
		\right\| 
		\leqslant \\ \leqslant
		\left\| 
	    		\sum\limits_{i=0}^{N-1}
			\int\limits_{\delta^{(i)}(t_0)}^{\delta^{(i+1)}(t_0)}
	    			\!\!\!\!\!\! 
			[\varphi(s, \xi) - \varphi(s, \xi_i)] \Delta{s}
		\right\| + \!\!\!\!\!
		\int\limits_{\delta^{(N)}(t_0)}^t
			\!\!\!\!\!
			\left\| \varphi(s, \xi(s)) \right\|\Delta{s}
		\leqslant
	\end{gathered}
    \]
    \[
	\begin{gathered}
		\leqslant
		\sum\limits_{i=0}^{N-1}
			\int\limits_{\delta^{(i)}(t_0)}^{\delta^{(i+1)}(t_0)}
	    			\!\!\!\!\!\! 
				\left\|\varphi(s, \xi) - \varphi(s, \xi_i)\right\|
				\Delta{s}
		+ 2 M \left(t - \delta^{(N)}(t_0)\right)
		\leqslant \\ \leqslant
		\sum\limits_{i=0}^{N-1}
			2\lambda 
	    		\!\!\!\!\!\! 
			\int\limits_{\delta^{(i)}(t_0)}^{\delta^{(i+1)}(t_0)}
	    		\!\!\!\!\!\! 
			\left\|\xi(s) - \xi_i\right\| \Delta{s}
		+ 2 M \left(\delta^{(N+1)}(t_0) - \delta^{(N)}(t_0)\right)
		\leqslant \\ \leqslant
		\sum\limits_{i=0}^{N-1}
			2\lambda \cdot \varepsilon M
	    		\left(\delta^{(i+1)}(t_0) - \delta^{(i)}(t_0)\right)
		+ 2MK
		\leqslant
	\end{gathered}
    \]
    \[
	\begin{gathered}
		\leqslant 
		2\lambda \cdot \varepsilon M
		\sum\limits_{i=0}^{N-1}
			\left(\delta^{(i+1)}(t_0) - \delta^{(i)}(t_0)\right)
		+ 2MK
		= \\ =
		2\lambda \cdot \varepsilon M
		\left(\delta^{(N)}(t_0) - t_0\right) + 2MK
		= \\ =
		2\lambda \cdot \varepsilon M \cdot \frac{L}{\varepsilon} + 2MK 
		= 2M\left(\lambda L + K\right).
	\end{gathered}
    \]

    Thus we have
    \begin{eqnarray*}
        \norm{x(t) - \xi(t)} & \leqslant & 
	\lambda\varepsilon\int\limits_{t_0}^t\norm{ x(s) -  \xi(s)}\Delta{s} +
               \varepsilon\norm{\int\limits_{t_0}^t
	       \varphi(s, \xi(s)) \Delta{s}} \leqslant \\
	& \leqslant & 
	\lambda\varepsilon\int\limits_{t_0}^t\norm{ x(s) -  \xi(s)}\Delta{s} + 
	\varepsilon \cdot 2M\left(\lambda L + K\right).
    \end{eqnarray*}
    
    Taking into account Gronwall's inequality and properties of the exponential
    function on time scale (\cite{BohnerPetersonAdvances}), we obtain as we did
    before
    \[
        \norm{x(t) - \xi(t)}  \leqslant 
		\varepsilon \cdot 
		2M\left(\lambda L + K\right) \cdot 
			e_{\lambda\varepsilon}(t, t_0) < 
		\varepsilon \cdot 
		2M\left(\lambda L + K\right) \cdot e^{\lambda L}, 
    \]
    that is,
    \[
	    \norm{x(t) - \xi(t)} < C\varepsilon,
    \]
    where $C = 2M\left(\lambda L + K\right) \cdot e^{\lambda L}$ 
    and this concludes the proof.
    \hfill $\quad \blacksquare$

It is clear that trivial time scales $\mathbb{R}$, $\mathbb{Z}$, and
$h\mathbb{Z}$ are periodic in shifts $\delta_{\pm}(T,t) = T \pm t$ for various
periods $T$.  Also, any periodic in shifts $\delta_{\pm}(T, t)$ function is
$\Delta$-periodic in such cases. Moreover, condition 3) of the last theorem is
trivially satisfied.  Thus proved theorem is the closest analogue of the averaging
theorem for ordinary differential equations with a periodic right-hand side.

At the same time to find a good example of periodic in shifts non-trivial time
scales appears to be a hard problem. Finding $\Delta$-periodic functions
defined on such time scales is a yet harder problem. For example, consider some
non-trivial time scale with a condensation point.  By definition,
a $\Delta$-periodic function has to compensate decreasing length of the
integration interval by increasing magnitude. Hence function needs to be
unbounded as time tends to condensation point and we cannot apply averaging
theorem.

\begin{example}
	Let $\mathbb{T} = \left\{t_n = 1 - \dfrac{1}{q^n}, n \in \mathbb{N}_0,
	q > 1\right\} \cup \left\{1\right\}$. This is a time scale with
	condensation point $t = 1$, forward jump operator $\sigma(t) = 
	\frac{q	- 1 + t}{q}$, and graininess $\mu(t) = \dfrac{q-1}{q}(1-t)$.
	Forward shift can be defined as follows:
	\[
		\delta_+(T, t) = \frac{q^T + t - 1}{q^T}.
	\]
	It is easy to compute $\delta^\Delta_+(T, t) = q^{-T}$.
	We found out a simple function $f(t) = \dfrac{1}{1-t}$ such that
	$f\left(\delta_+(T, t)\right)\delta^\Delta_+(T, t) = f(t)$, i.~e. 
	the function $f(t)$ is $\Delta$-periodic in shifts. However $f(t)$ is
	unbounded above as $t \to 1$.
\end{example}

Analyzing the example, we found one more possibility to obtain a more accurate
estimate for proximity between solutions of the original and averaged systems.
\begin{definition}
	Let $\mathbb{T}$ be a periodic in shift $\delta_+(P, t)$ time scale with a period $P$.
	A function $f(t)$ is called geometric $\Delta$-quasiperiodic function with period $T > P$
	and factor $\gamma$ if the following condition holds:
	\begin{equation}
		\label{quasiperiodic}
		f\left(\delta_+(T, t)\right)\delta^\Delta_+(T, t) = \gamma f(t).
	\end{equation}
\end{definition}

Using substitution rule \eqref{substitution-rule} we can easily prove 
the important property of geometric $\Delta$-quasiperiodic function.

\begin{lemma}
	\label{quasiperiodic-substitution-rule}
	Let $\mathbb{T}$ be a time scale that is periodic in shift $\delta_+$
	with the period $P$ and $f$ a geometric $\Delta$-quasiperiodic in shift $\delta_+$ 
	with the period $T \in [P, \infty)_{\mathbb{T}^*}$.
	Suppose that $f \in C_{rd}(\mathbb{T})$, then
	\[
		\int\limits_{t_0}^{t} f(s) \Delta{s} = 
		\gamma \int\limits_{\delta^{T}_+(t_0)}^{\delta^{T}_+(t)} f(s) \Delta{s}.
	\]
\end{lemma}
\begin{proof}
Substituting $\nu(s) = \delta_+(T, s)$ and $g(s) = f\left(\delta_+(T, t)\right)$ in 
\eqref{substitution-rule} and taking \eqref{quasiperiodic} into account we obtain 
the statement of lemma by direct calculation.
\end{proof}

Now suppose $X(t, x)$ in \eqref{mainsys-periodic} is geometric $\Delta$-quasiperiodic 
with period $T$ and factor $\gamma$ for any fixed $x$. Consider dynamic system
\begin{equation}
    \label{avgsys-quasiperiodic}
    \xi^{\Delta} = \varepsilon \widehat{X}(t, \xi), \quad \xi(t_0) = x_0,
\end{equation}
where
\begin{equation}
    \label{avglimit-quasiperiodic}
    \begin{aligned}
    \widehat{X}(t, x) = 
    	\Bigg\{
	\widehat{X}_i(x) 
		&= \frac{\gamma^i}{\delta^{(i+1)}(t_0) - \delta^{(i)}(t_0)} 
			\int\limits^{\delta_+(T, t_0)}_{t_0} X(t,x) \Delta{t}, \\
		& \delta^{(i)}(t_0) \leqslant t < \delta^{(i+1)}(t_0), 
			\quad i = 0, 1, 2, \dots 
	\Bigg\}.
    \end{aligned}
\end{equation}

We can prove now that there exists proximity between 
solutions of systems \eqref{mainsys-periodic} and \eqref{avgsys-quasiperiodic} when
$X(t, x)$ is a geometric $\Delta$-quasiperiodic function.

\begin{theorem}
\label{theorem-quasiperiodic-averaging}
Suppose the conditions 2)--4) of Theorem \ref{theorem-periodic-averaging} hold
in $Q$, and besides this, every component of vector--function $X(t, x)$ is
geometric $\Delta$-quasiperiodic function with period $T$ and factor $\gamma$ for any
fixed $x$. 

Then for any $L > 0$ there exists $\varepsilon_0\left(L\right) > 0$ such that
for $0 < \varepsilon <~\varepsilon_0$ and $t \in \left[t_0, t_0 + L\varepsilon^{-1}\right] \cap \mathbb{T}$ 
the following estimate holds:
\begin{equation}
	\label{solution-quasiproximity}
	\norm{x(t) - \xi(t)} \leqslant C\varepsilon,
\end{equation}
where $x(t)$ and $\xi(t)$ denote solutions of the Cauchy problems
\eqref{mainsys-periodic} and \eqref{avgsys-quasiperiodic} respectively. 

\end{theorem}
\begin{proof}
From quasiperiodical properties of the function $X(t, x)$,
it follows easily that 
\[
    \int\limits_{\delta^{(i)}(t_0)}^{\delta^{(i+1)}(t_0)} 
	\!\!\!\!\!\!\varphi(s, \xi_i) \Delta{s} = 0, \qquad i = 0, 1, \dots,
\]
where $\varphi(t, \xi) = X(t, \xi(s)) - \widehat{X}(t, \xi(s))$ and
$\xi_i = \xi\left(\delta^{(i)}(t_0)\right) = const$. 
Thus the argumentation of previous proof can be repeated almost literally.
For brevity, we omit the details.
\end{proof}

\begin{example}
	Let us use time scale from previous example. Consider dynamic system
	\[
		x^\Delta = \varepsilon (-1)^{-\frac{\ln(1-t)}{\ln{q}}} x, 
			\quad x(0) = 1, \quad t \in \mathbb{T}, 
	\]
	that is, $X(t, x) = (-1)^{-\frac{\ln(1-t)}{\ln{q}}} x$.
	We get 
	\[
	\begin{aligned}
		X\left(\delta_+(T, t), x\right)\delta^\Delta_+(T, t) &= 
		x \cdot (-1)^{-\frac{\ln\left(1 - \frac{q^T + t - 1}{q^T}\right)}{\ln{q}}} \cdot \frac{1}{q^T} =\\ 
		&= x \cdot (-1)^{-\frac{\ln(1 - t) - \ln{q^T}}{\ln{q}}} \cdot \frac{1}{q^T} = \\
		&= x \cdot (-1)^{-\frac{\ln(1 - t)}{\ln{q}}} \cdot (-1)^T \cdot \frac{1}{q^T} = \\
		&= X(t, x) \cdot (-1)^T \cdot \frac{1}{q^T}.
	\end{aligned}
	\]
	This implies that $X(t, x)$ is geometric $\Delta$-quasiperiodic with
	period $T = 2$ and factor $\gamma = q^{-T} = q^{-2}$.

	Further, $\delta^{(i+1)}(0) - \delta^{(i)}(0) = \dfrac{q^T - 1}{q^{T(i+1)}} \leqslant \dfrac{q^T - 1}{q^T}$.
	Thus we have
	\[
		\begin{aligned}
		\widehat{X}_i(x) &= 
			\frac{\left(q^{-T}\right)^i}{\delta^{(i+1)}(0) - \delta^{(i)}(0)} 
			\int\limits^{\delta_+(T, 0)}_{0} X(t,x) \Delta{t} =  \\
			&= \frac{q^2}{q^2 - 1} \int\limits^{1 - \frac{1}{q^2}}_{0} X(t,x) \Delta{t} =  \\
			&= \frac{q^2}{q^2 - 1} \cdot x \cdot \frac{q-1}{q^2} = x \cdot \frac{1}{q+1}.
		\end{aligned}
	\]
	Hence we have two systems on the same time scale:
	\[
		\left\{\begin{aligned}
			&x^\Delta = \varepsilon \cdot (-1)^{-\frac{\ln(1-t)}{\ln{q}}} x, \\
			&x(0) = 1, 
		\end{aligned}\right.
		\qquad 
		\text{and}
		\qquad
		\left\{\begin{aligned}
			&\xi^\Delta = \varepsilon \cdot \dfrac{\xi}{q + 1}, \\
			&\xi(0) = 1. 
		\end{aligned}\right.
	\]
	It is not too hard to find exact solution of the linear equation 
	\[
		y^\Delta = p y, \quad y(0) = y_0, \quad t \in \mathbb{T}.
	\]
	Indeed, all $t \neq 1$ are isolated points and thus $y(\sigma(t)) = y(t) + \mu(t)y^\Delta(t)$.
	Starting from $t = 0$ we get $y\left(\sigma^k(0)\right) = 
	y_0 \prod\limits_{i=0}^{k-1}\left[1 + p\mu\left(\sigma^i(0)\right)\right]$.
	This yields that
	\[
		y(t) = y_0 \prod\limits_{i=0}^{k-1}\left[1 + \dfrac{p(q-1)}{q^{i+1}}\right], 
			\quad k = -\frac{\ln(1-t)}{\ln{q}}, \quad t \neq 1.
	\]
	Actually $y(t) = e_p(t, 0)$, i. e. exponential function on time scale $\mathbb{T}$.

	In the same way, we obtain exact solutions of original and averaged systems:
	\[
		\begin{gathered}
			x(t) = \prod\limits_{i=0}^{k-1}\left[1 + \varepsilon \cdot \dfrac{(-1)^i (q - 1)}{q^{i+1}}\right], \\
			\xi(t) = \prod\limits_{i=0}^{k-1}\left[1 + \varepsilon \cdot \dfrac{q - 1}{q^{i+1}(q+1)}\right]. \\
		\end{gathered}
	\]
	It seems to be impossible to find a precise analytical estimate of difference
	$\left|x(t) - \xi(t)\right|$ in terms of $\varepsilon$. Instead we
	conducted numerical modelling and found empirical dependence between
	proximity of solutions and small parameter $\varepsilon$. The results of
	modelling are presented in Figure 1.
	
	\begin{figure}[h!]
		\begin{center}
			\begin{minipage}{0.48\textwidth}
				\centering
				\includegraphics[scale=0.25]{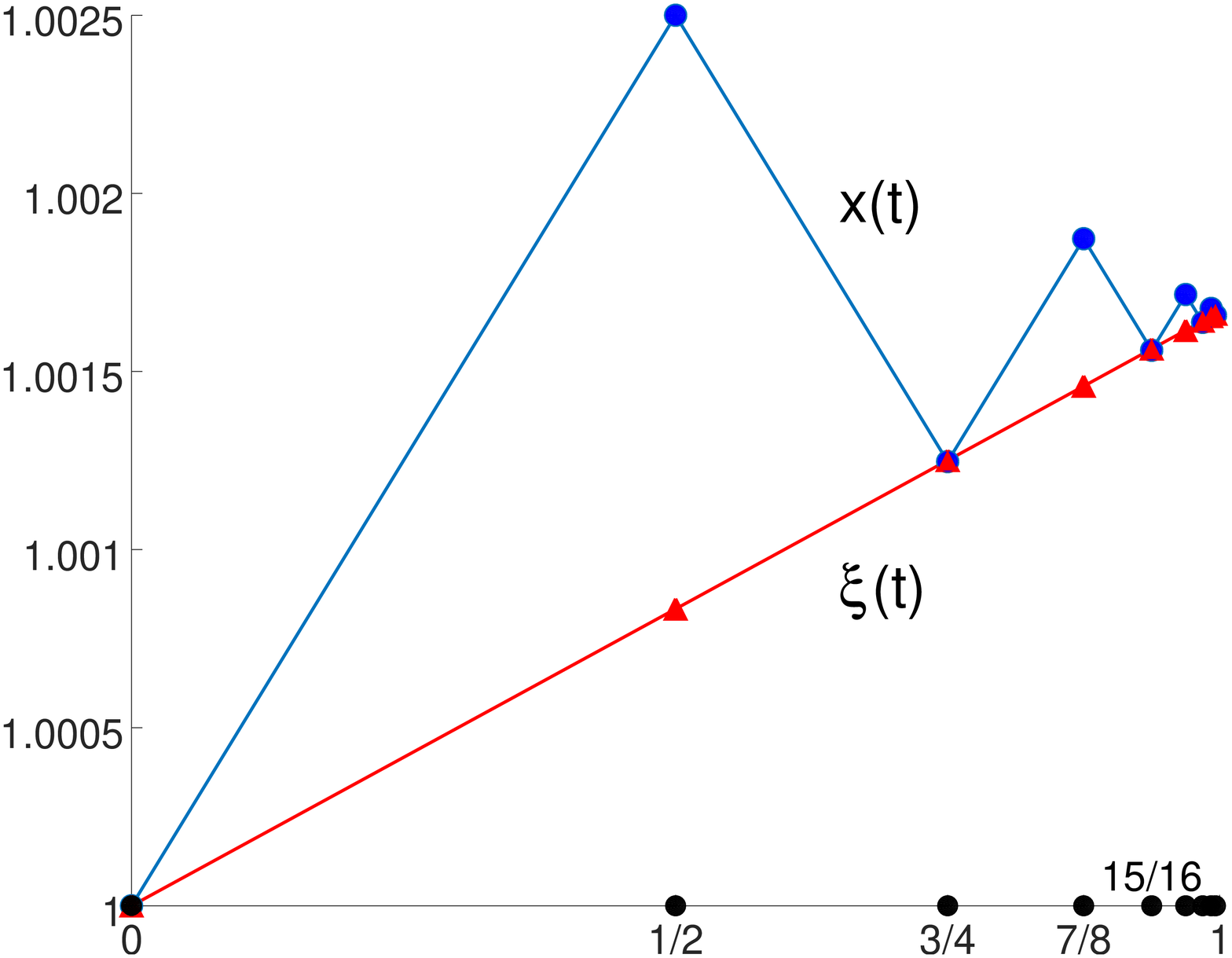} \\
				{\footnotesize a) Solutions of original and averaged systems,
					$\varepsilon = 0.005$, $q = 2$}
			\end{minipage}
			\hspace{3 mm}
			\begin{minipage}{0.48\textwidth}
				\centering
				\includegraphics[scale=0.25]{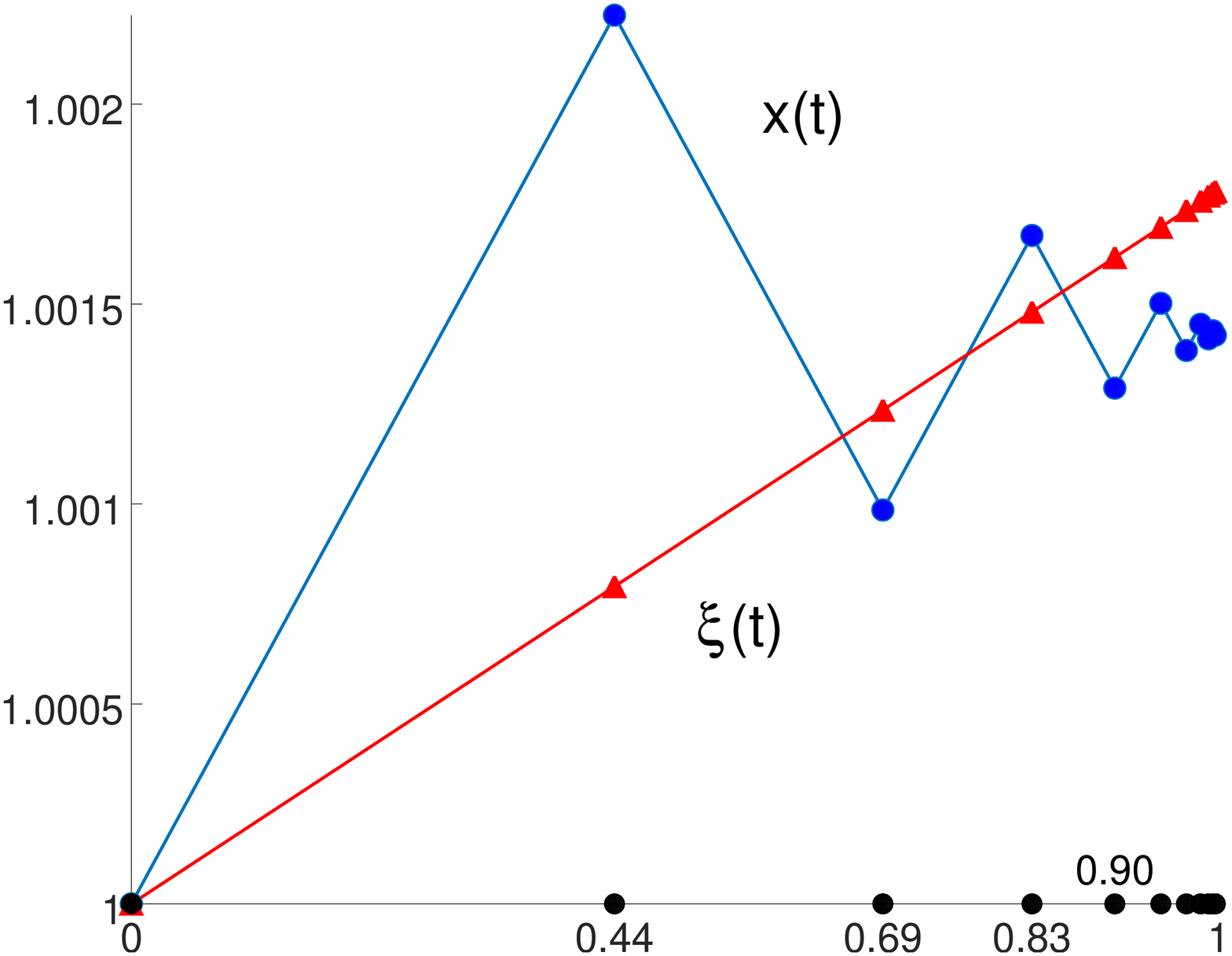} \\
				{\footnotesize a) Solutions of original and averaged systems,
					$\varepsilon = 0.005$, $q = 1.8$}
			\end{minipage}
		\end{center}
		\begin{center}
			\begin{minipage}{0.48\textwidth}
				\centering
				\includegraphics[scale=0.25]{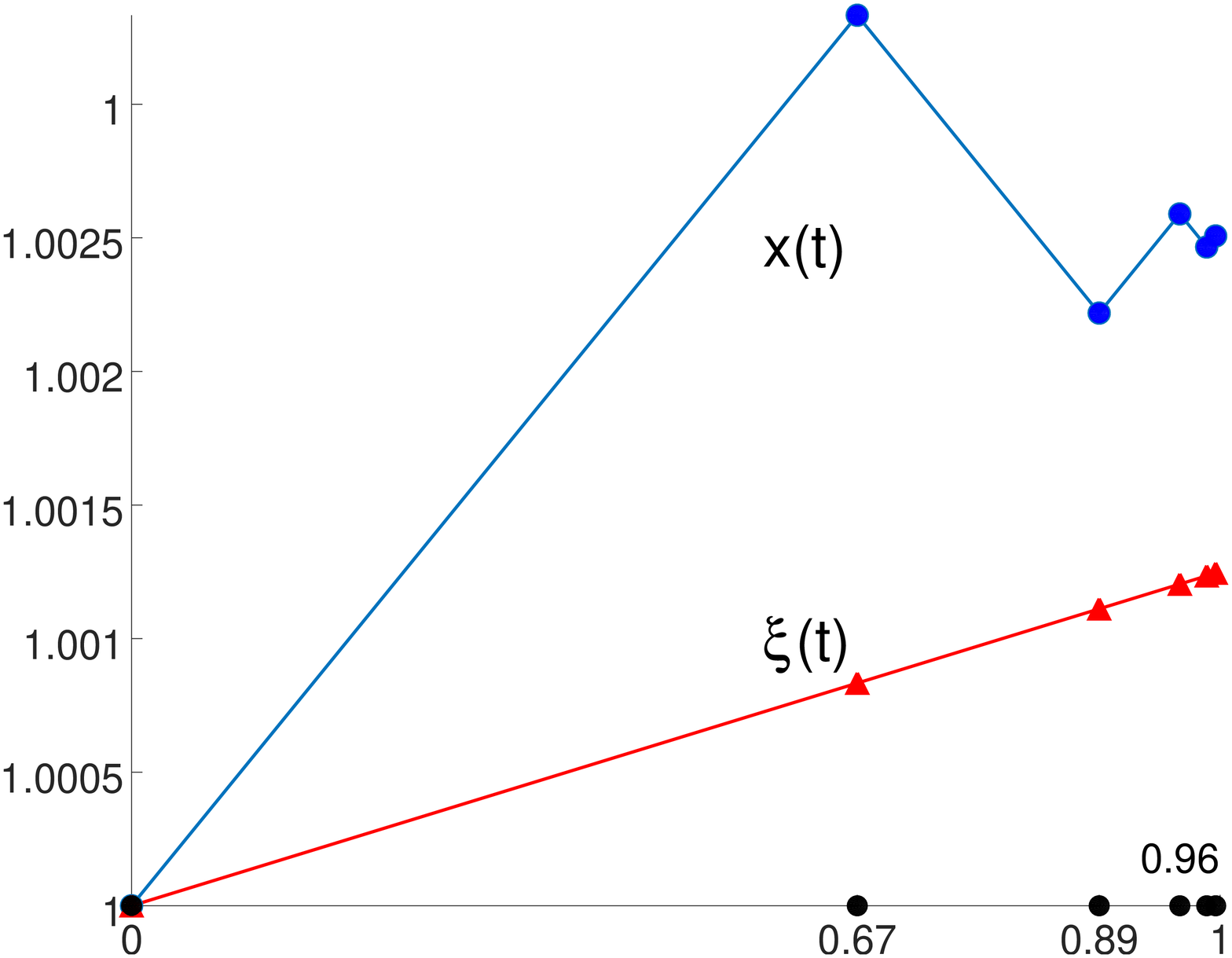} \\
				{\footnotesize a) Solutions of original and averaged systems,
					$\varepsilon = 0.005$, $q = 3$}
			\end{minipage}
			\hspace{3 mm}
			\begin{minipage}{0.48\textwidth}
				\centering
				\includegraphics[scale=0.25]{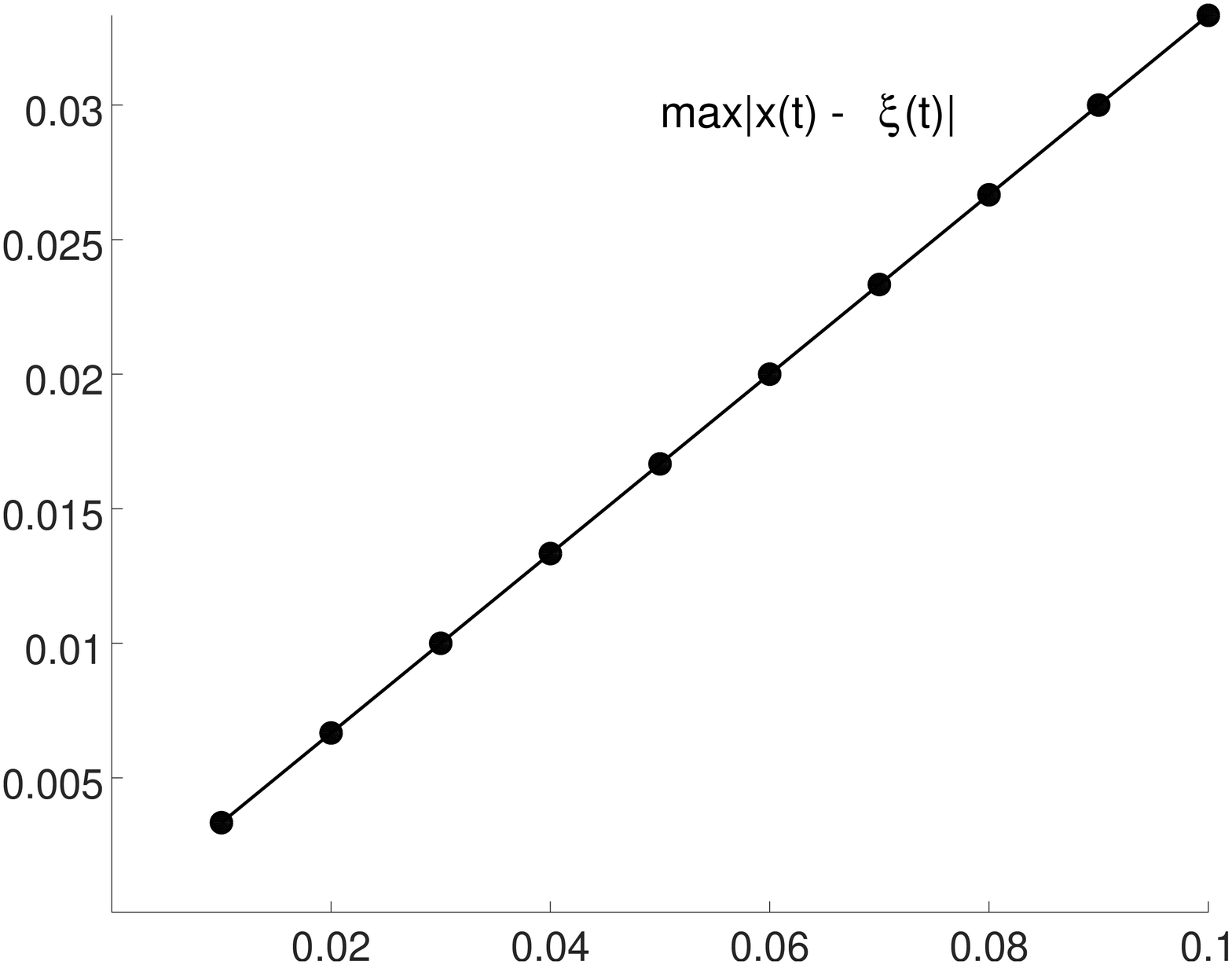} \\
				{\footnotesize d) Absolute difference between solutons in regard to $\varepsilon$}
			\end{minipage}
		\end{center}
		\caption{Numerical modelling of averaging method for quasiperiodic system on
			time scale $\mathbb{T} = \left\{t_n = 1 - \dfrac{1}{q^n}, n \in \mathbb{N}_0,
			q > 1\right\} \cup \left\{1\right\}$.
		}
	\end{figure}
\end{example}

\section{Conclusion.}
The aim of this paper is to develop our previous results for the averaging
method on time scales. Following \cite{Adivar2013} we considered
$\Delta$--periodic systems and obtained a more accurate estimate for proximity
between solutions of original and averaged systems. Moreover, the same result
was obtained for dynamic systems with a quasiperiodic right-hand side, which
are introduced for the first time.  To illustrate the application of the
averaging theorem for such kind of system we considered an example and
conducted numerical modelling. Obtained results can be used to improve
previously developed numer\-ical\-ly--asymptotic method of solution for optimal
control problems on time scales.

\bibliographystyle{elsarticle-harv} 
\bibliography{references}

\end{document}